\documentstyle[11pt]{article}
\begin{document}
\newcommand{\be}{\begin{equation}}
\newcommand{\ee}{\end{equation}}
\newcommand{\bea}{\begin{eqnarray}}
\newcommand{\eea}{\end{eqnarray}}
\newcommand{\beaa}{\begin{eqnarray*}}
\newcommand{\eeaa}{\end{eqnarray*}}
\newcommand{\qd}{\quad}
\newcommand{\qqd}{\qquad}
\newcommand{\npb}{\nopagebreak[1]}
\newcommand{\nn}{\nonumber}
\title{\bf Generalized Hirota bilinear identity\\
and integrable q-difference \\and lattice hierarchies}
\author{L.V.Bogdanov\thanks{e-mail leonid@landau.ac.ru}
\\   Landau Institute for
Theoretical Physics,
IINS,\\
Kosygina str. 2, Moscow v-334, GSP-1 117940, Russia
}
\date{}
\maketitle
\begin{abstract}
Hirota bilinear identity for Cauchy-Baker-Akhieser (CBA) kernel
is introduced as a basic tool to construct
integrable hierarchies containing lattice and q-difference times.
Determinant formula for the action of meromorphic function
on CBA kernel is derived. This formula gives opportunity
to construct generic solutions for integrable lattice
equations.
\end{abstract}
\section{Introduction.}
This paper is a sequel of the work \cite{Bogdanov}.
We introduce here generalized Hirota bilinear identity and develop
a unified approach to continuous, lattice and
q-difference variables in integrable hierarchies.

Generalized Hirota bilinear identity uses the function of
two complex variables $\psi(\lambda,\mu)$, which is an
analogue of Cauchy-Baker-Akhieser kernel
on the Riemann
surface  \cite{Orlov} in Segal-Wilson type Grassmannian context
\cite{Segal}. In frame of $\bar{\partial}$-dressing method
this function is connected with the solution of nonlocal
$\bar{\partial}$-problem normalized by $(\lambda-\mu)^{-1}$
(see \cite{Manakov}, \cite{BM}).
The action of meromorphic function on the Grassmannian
can be explicitly found in terms of the function $\psi(\lambda,\mu)$
in the
form of elegant determinant formula, which has close ties
with the Miwa's formula for the $\tau$-function in the
model of free fermions \cite{Miwa} and Wronskian formulae for the
composition of Darboux transformations.

Lattice versions and q-deformations of KP hierarchy
and N-waves system are constructed. The determinant
formula gives generic solution (i.e. solution corresponding
to arbitrary Grassmannian point) in explicit form.
The same is valid for q-difference equations, but in this case,
however,
the determinants of infinite-dimensional matrices appear in the
formula.
\section{Generalized Hirota identity.}
Due to the limited volume of this publication,
we will take quite a formal start, introducing from
the beginning generalized Hirota bilinear identity as a basic tool
for the following consideration
\be
\int_{\partial G} \chi(\lambda,\nu;g_1)g_1(\nu)g_2^{-1}(\nu)
\chi(\nu,\mu;g_2)d\nu=0 ,
\label{HIROTA}
\ee
Here $\chi(\lambda,\mu;g)$ is a function of two complex variables
$\lambda ,\mu\in G$ and a functional of the group element $g$
defining the dynamics (will be specified later), $G$ is some set
of domains of the complex plane, the integration goes over the
boundary of $G$. By definition, the function $\chi(\lambda,\mu)$
possesses the following analytical properties: as $\lambda
\rightarrow\mu$, $\chi\rightarrow (\lambda-\mu)^{-1}$ and
$\chi(\lambda,\mu)$ is analytic function of
two variables $\lambda,\mu$ for $\lambda\neq\mu$.

In  another form, more similar to standard Hirota bilinear identity,
the identity (\ref{HIROTA}) can be written as
\be
\int_{\partial G} \psi(\nu,\lambda,g_1)
\psi(\lambda,\mu,g_2)d\lambda=0 ,
\ee
where
$$
\psi(\lambda,\mu,g)=g^{-1}(\lambda)\chi(\lambda,\mu,g)g(\mu)
$$
The place of the function $\chi$ in the Segal-Wilson Grassmannian
approach is outlined in \cite{Bogdanov}. In frame of $\bar{\partial}$
-dressing method the analogue of this function is well-known
(see \cite{BM}), in algebro-geometric technique the function
$\psi(\lambda,\mu)$ corresponds to Cauchy-Baker-Akhieser kernel
on the Riemann surface (see \cite{Orlov}). So the function
$\psi(\lambda,\mu)$ introduced here could be called CBA kernel
for Segal-Wilson type Grassmannian. Let us consider two
linear spaces $W(g)$ and $W'(g)$ defined by the function $\chi(\lambda,\mu)$
via equations connected with the identity (\ref{HIROTA})
\bea
\int_{\partial G} \chi(\lambda,\nu;g)
f(\nu;g)d\nu=0,
\label{W}\\
\int_{\partial G}h(\nu;g)\chi(\nu,\mu;g)
d\nu=0 ,
\label{W'}
\eea
here
$f(\lambda)\in W$, $h(\lambda)\in W'$; $f(\lambda)$,
$g(\lambda)$ are defined in $\bar G$.
It is easy to check using the analytical properties
of the function $\chi(\lambda,\mu)$ that these spaces
possess the following properties

1. $W,W'$ contains a meromorphic function with arbitrary given
divisor of poles in $G$ (completeness)

2. $W,W'$ is transversal to the space of holomorphic functions
in $G$ (transversality).

So the function $\chi(\lambda ,\mu)$ defines simultaneously
a point and a dual point of `analytical Grassmannian',
where by `analytical Grassmannian' we understand the set of linear
spaces $W$ of
functions in $\bar G$ possessing the properties 1,2. The function
$\chi(\lambda,\mu)$ plays a role of the basic function
in both of the spaces $W,W'$ (for $W$ with respect to the
first variable, for $W'$ with respect to the second variable).
It follows from the definition of linear spaces $W,W'$ that
\bea
f(\lambda)&=&2\pi {\rm i}\int\!\!\!\int_G \left({\partial\over \partial\bar\nu}
f(\nu)\right)\chi(\lambda,\nu)d\nu\wedge d\bar{\nu},
\nn\\
g(\mu)&=&-2\pi {\rm i}\int\!\!\!\int_G \left({\partial\over \partial\bar\nu}
g(\nu)\right)\chi(\nu,\mu)d\nu\wedge d\bar{\nu},
\label{basis}
\eea
these formulae in some sense provide an expansion of the
functions $f,g$ in terms of the basic function $\chi(\lambda,\mu)$.

The dynamics of linear spaces $W,W'$ looks very
simple
\be
W(g)=gW_0;\quad W'(g)=g^{-1}W'_0,
\label{dynamics}
\ee
here $W_0=W(g=1)$, $W'_0=W'(g=1)$ (the formulae (\ref{dynamics})
follow from identity (\ref{HIROTA}) and the formulae (\ref{basis})).

\section{Determinant formula for action of rational
$g(\lambda)$ on the CBA kernel}
The action of arbitrary meromorphic function $g(\lambda)$,
having in $G$ equal number of zeroes and poles
(taking multiplicity into account), on the
function $\chi(\lambda,\mu)$ can be found in explicit form.
This result was obtained in our work \cite{Bogdanov},
using the properties of the space $W$.
In the present work we would like to start from
Hirota bilinear identity (\ref{HIROTA}) and introduce a
compact determinant formula for this action.

Let $g(\lambda)$ be some meromorphic
function in $G$ having simple poles at the set of points $z_i$
and simple zeroes at the set of point $\hat z_i$,
$0<i<N$, and
let the basic function  $\chi_0(\lambda,\mu)$
be defined at the initial point (i.e. for $g(\lambda)=1$).
The problem is to find the solution of identity (\ref{HIROTA})
\be
\int_{\partial G} \chi_0(\lambda,\nu)g^{-1}(\nu)
\chi(\nu,\mu;g)d\nu=0 ,
\label{HIROTAA}
\ee
having analytical properties specified in the definition
(see (\ref{HIROTA})). The answer is quite simple
\be
\chi(g)=
g(\lambda)g^{-1}(\mu)
{\Delta(\lambda,z_1,...,z_N;\mu,\hat z_1,...,\hat z_N)
(\chi_0(\lambda,\mu))
\over \Delta(z_1,...,z_N;\hat z_1,...,\hat z_N)
(\chi_0(\lambda,\mu))},
\label{Det}
\ee
where
\be
\Delta(z_1,...,z_N;\hat z_1,...,\hat
z_N)(f(\lambda,\mu))=\det(f_{ij}):=
\det(f(z_i, \hat z_j)).
\ee
The function $\chi(\lambda,\mu,g)$
defined by the formula (\ref{Det}) satisfies the equation
(\ref{HIROTAA}) and possesses necessary analytic properties.
The most simple way to derive this formula is to use the
consequence of the Hirota bilinear identity (\ref{HIROTA})
\be
\int\!\!\!\int_{\partial G\times\partial G} \chi_0(\lambda,\nu)g^{-1}(\nu)
\chi(\nu,\eta;g)g(\eta)\chi_0(\eta,\mu)d\eta d\nu=0.
\ee

For the arbitrary meromorphic function $g(\lambda)$, having
poles with the multiplicity $n_i$ at the set of points $z_i$
and zeroes with the multiplicity $\hat n_i$
at the set of point $\hat z_i$
\be
\chi(g)=g(\lambda)g^{-1}(\mu)
{\Delta((\lambda,1),(z_1,n_1),...,(z_N,n_N);(\mu,1),
(\hat z_1,\hat n_1),...,(\hat z_N,\hat n_N))
(\chi_0(\lambda,\mu))
\over \Delta((z_1,n_1),...,(z_N,n_N);(\hat z_1,\hat
n_1),...,(\hat z_N,\hat n_N))
(\chi_0(\lambda,\mu))},
\label{Det1}
\ee
where
\be
\Delta((z_1,n_1)...,(z_N,n_N);(\hat z_1,\hat n_1),...,(\hat z_N,\hat n_N))
(\chi(\lambda,\mu))={\rm det}(\chi_{IJ}),
\ee
the matrix $\chi_{IJ}$
consists of the elements
$$\left({\partial^{p_i-1}\over\partial\lambda^{p_i-1}}
{\partial^{\hat p_j-1}\over\partial\mu^{\hat
p_j-1}}\chi\right)(z_i,\hat
z_j)=\chi_{IJ},$$
here $0<p_i\leq n_i,\:0<\hat p_j\leq\hat n_j$, $I=\sum_{k=0}^{i-1}n_i+
p_i$, $J=\sum_{k=0}^{i-1}\hat n_j+\hat p_j$.

The function $\chi(\lambda,\mu)$ is connected with the $\tau$-function
by the formula
\be
\chi(\nu,\mu)={\tau\left(g\times \left({\lambda-\nu\over\lambda-\mu}
\right)\right)\over\tau(g)(\nu-\mu)}
\label{TAUFORM}
\ee
We do not give the direct proof of the formula (\ref{TAUFORM})
, though it seems more
or less evident if $G$ is a unit circle.
But for not 1-connected $G$ or $G$ consisting of several disconnected
domains the interpretation of the right part is not so trivial,
it goes beyond the standard Segal-Wilson approach.
The formula (\ref{TAUFORM}) gives an opportunity to obtain
Miwa's formula \cite{Miwa} as a special case of the formula
({\ref{Det1}).

\section{Introduction of lattice and q-difference variables
to integrable hierarchies}

A dependence of the function
$\chi(\lambda,\mu)$ on dynamical variables
is hidden in the function $g(\lambda)$.
Usually these variables are continuous space and
time variables, but it is possible also to introduce discrete (lattice)
an q-difference variables into identity (\ref{HIROTA}).
We will consider the
following functions $g(\lambda)$
\begin{eqnarray}
g^{-1}_i&=&\exp(K_i x_i)\label{c}
;\quad
{\partial\over\partial x_i}g^{-1}=K_i g^{-1},\\
g^{-1}_i&=&(1+l_i K_i)^{n_i}\label{d};\:
\Delta_i g^{-1}={g^{-1}(n_i+1)-g^{-1}(n_i)\over
l_i}=K_i g^{-1},
\\
g^{-1}_i&=&{\rm e}_q (K_i y_i);\label{q}\quad
\delta_i^q g^{-1}={g^{-1}(qy_i)-g^{-1}(y_i)\over
(q-1)y_i}=K_i g^{-1}.
\end{eqnarray}
Here $K_i(\lambda)$ are rational functions. The function
(\ref{c}) introduces a dependence on continuous variable
$x_i$, the function (\ref{d}) -- on discrete variable $n_i$
(the lattice parameter $l_i$ may depend on $n_i$, so it is possible
to consider lattices with the changing step of the lattice)
and the function (\ref{q}) defines a dependence of $\chi(\lambda)$
on the variable $y_i$ (we will call it a q-difference variable).
The function ${\rm e}_q(y)$ has a representation
\be
{\rm e}_q(y)=(\prod_{n=1}^{\infty}(1+q^ny(q-1)))^{-1},
\label{qexp}
\ee
here we suggest that $|q|<1$.

To introduce a dependence on several variables (may be of different
type), one should consider a product of corresponding functions $g(\lambda)$
(all of them commute).

The structure of the functions (\ref{c}-\ref{q}) imposes
some limitations on the choice of $G$ in order to apply
formula (\ref{Det1}) to construct solutions of corresponding
lattice and q-difference equations. Typically $G$ should contain
all zeroes and poles of the considered class of functions $g$.
In the case of function (\ref{d}) it is just a finite set
of neighborhoods of points, in the case (\ref{q}), $q\in{\bf R}$
$G$ should contain neighborhoods of some curves in the complex
plane, where group elements (\ref{q}) have zeroes.
Every function $K_i(\lambda)$
may be defined on its own copy of the complex plane; then
it needs a definition on another's copy (we will use zero value
on another's copy). However, the choice of $G$ is not so important
for construction of equations.
\section{Construction of equations}
Equations in the right part of
(\ref{c}-\ref{q}) and the boundary condition $g(0)=1$
characterize the corresponding functions (and give a definition
of ${\rm e_q(y)}$).
These equations
play a crucial role in the algebraic scheme of constructing
integrable equations.

This scheme is based on the assumption of solvability
of the problem (\ref{HIROTA}) for some class
of functions $g(\lambda)$, on transversality property
and on the existence of special
operators, which transform $W$ into itself.

Indeed, the relation (\ref{W}) implies that
if $\chi({\bf x},{\bf n},{\bf y},\lambda)\in W({\bf x},{\bf n},
{\bf y})$,
then the functions
\bea
D_{i}^c \chi&=&\partial / \partial x_{i}\chi +\chi K_i(\lambda) \nn\\
D_{i}^d \chi&=&\Delta_{i}\chi +T_i^d \chi K_i(\lambda) \nn\\
D_{i}^q \chi&=&\delta^q_{i}\chi +T_i^q \chi K_i(\lambda)
\label{D}
\eea
also belong to $W$, where $Tf(n)=f(n+1)$,
$T^q f(y)=f(qy)$. We can multiply the
solution from the left
by the arbitrary matrix function
of additional variables,
$u({\bf x},{\bf n},{\bf y})\chi
\in W$. So the operators (\ref{D}) are the generators
of Zakharov-Manakov ring of operators, that transform
$W$ into itself.

Combining this property
with the transversality property,
one obtains the differential relations
between the coefficients of expansion of functions
$\chi({\bf x},{\bf n}, {\bf y},\lambda)$ into powers of
$(\lambda-\lambda_{p})$ at the poles of
$K_{i}(\lambda)$ \cite{BM}.

The derivation of equations in this case is completely
analogous to the continuous case \cite{BM}.
\subsection{Lattice and q-difference KP hierarchy}
The following derivation will be conducted for q-difference case,
to get the difference case you should just change
$\delta^q_i$ for $\Delta_i$ and $T_i$ for $T^q_i$.

The KP hierarchy corresponds to the special choice
of the functions $K_i(\lambda)$: $K_i(\lambda)=(l_i)^{-1}
\lambda^{-i}$.
The transversality of
the space $W$ implies linear equations (below we transformed
operators $D^q$ to $\delta^q$ by the substitution
$\psi(\lambda)=g(\lambda)\chi(\lambda,0)$)
\be
(\delta^q_i-{\delta^q_1}^i)\psi(\lambda)=
\sum_{k=0}^{i-1}u_k ({\bf y}){\delta^q_1}^k \psi(\lambda).
\ee
The q-difference KP hierarchy is the set of compatibility conditions
for these linear equations.

\subsection{q-difference N-waves equations}
An arbitrary rational function $K_i(\lambda)$ may be
treated as degenerate case of the function having simple poles.
So generic equations in the hierarchy (\ref{q})
correspond to the operators $D_i$ with simple and
distinct poles. These equations play a fundamental role,
they can be constructed explicitly (see \cite{BM} for the
continuous case).

We use the functions $K_{i}(\lambda)$ :
\begin{equation}
K_{i}(\lambda)=\sum_{\alpha=1}^{n_{i}}\frac{a_{i}^{\alpha}}
{\lambda - \lambda_{i}^{\alpha}},
\end{equation}
where $a_{i}^{\alpha},\lambda_{i}^{\alpha} \in {\bf C}$,
$1 \leq \alpha \leq n_{i}$, $\lambda_{i}^{\alpha} \neq
\lambda_{j}^{\beta}$.
Then due to transversality of the space $W$
the function $\chi(\lambda,\mu)$ satisfy the relations
\begin{eqnarray}
D_{i}^q \chi(\lambda,\mu,{\bf y}) -
K_{i}(\mu)\chi(\lambda,\mu,{\bf y})-
(T_i^q\chi(\lambda_i^{\alpha},\mu,{\bf y}))
\chi(\lambda,\lambda_
i^{\alpha}, {\bf y})=0,
\label{N}
\end{eqnarray}
summation over $\alpha$ is understood. The other way to derive
this equation is to use Hirota bilinear identity
(\ref{HIROTA}).

The substitution of the values
$\lambda =\lambda_{k}^{\gamma}$, $\mu=\mu_{j}^{\beta}$,
$i\ne j \ne k\ne i$  to (\ref{N}) yields the equation
\begin{eqnarray}
\delta_i^q \chi_{jk}^{\beta
\gamma}+ K_{i}(\lambda_{k}^{\gamma})T_i^q\chi_{jk}^
{\beta \gamma}
-K_{i}(\lambda_{j}^{\beta})\chi_{jk}^{\beta \gamma}-
(T_i^q\chi_{ji}^{\beta \alpha})a_{i}^{\alpha} \chi_{ik}^{\alpha
\gamma}=0,
\label{NW}\\
\chi_{jk}^{\beta \gamma}=\chi(\lambda_k^{\gamma},
\lambda_j^{\beta},{\bf y})\nn
\end{eqnarray}
summation over $\alpha$ is understood.
If different permutations $ijk$ and substitutions of
the indices $\beta ,\gamma$ are taken into account, (\ref{NW})
is a closed set of equations for the functions
$\chi_{ji}^{\beta \alpha}({\bf n})$ for chosen $ijk$.

Expression (\ref{Det1}) gives the solutions
for the system of equations (\ref{NW}) and also for the
lattice KP hierarchy
at the arbitrary point of the lattice. The function
$\chi_0(\lambda,\mu)$ defines the point of the Grassmannian
and plays the role of spectral data.

The function $g(\lambda)$ (\ref{q}), taking into account
the representation (\ref{qexp}), has infinite number of zeroes
and poles. So application of the formula (\ref{Det1}) to
q-difference case requires the use of determinants of
infinite-dimensional matrices.

We would like also to consider a simple example in which
every function $K_i(\lambda)$ is defined on its own copy
of the complex plane.
Let's consider the lattice case.
We take a set of unit disks
${\bf D}_i$
as $G$.
The functions $K_i(\lambda)$ are chosen in the form
\bea
K_i(\lambda)={1\over \lambda},\quad \lambda\in{\bf D}_i;\nn\\
K_i(\lambda)=0,\quad \lambda\notin {\bf D}_i.\nn
\eea

We could use Zakharov-Manakov ring of operators to derive
the equations, but it is easy in this case to obtain them
directly from the identity (\ref{HIROTA}). Performing the integration
in the formula
\be
\int_{\partial G} \chi(\lambda,\nu,{\bf n})(1+K_i(\nu)l_i)
T_i\chi(\nu,\mu, {\bf n})d\mu=0 ,
\ee
one obtains
\bea
(1+K_i(\lambda)l_i)T_i\chi(\lambda,\mu, {\bf n})-
\chi(\lambda,\mu,{\bf n})(1+K_i(\mu)l_i)=\nn\\
\chi(\lambda,0_i,{\bf n})l_i
T_i\chi(0_i,\mu, {\bf n})=0
\label{N1}
\eea
where $0_i$ is a zero point of the corresponding region.
Taking the equation (\ref{N1}) at $\lambda=0_j,\mu=0_k$,
$i\ne j\ne k\ne i$
we get
\be
\Delta_i \chi_{jk}=(T_i\chi_{ji})
\chi_{ik}
\label{N2}
\ee
where $\chi_{ik}({\bf n})=\chi(0_k,0_j,{\bf n})$
If different permutations $ijk$ are taken into account, (\ref{N2})
is a closed set of equations for the functions
$\chi_{ji}({\bf n})$.
\subsection*{Acknowledgments}
This work was supported in part
by INTAS (International
Association for the promotion of cooperation with scientists
from independent states of the former Soviet Union),
Soros Foundation (ISF) (grant MLY 000)
and Russian Foundation for Fundamental Studies.

\end{document}